\documentclass[12pt,twoside]{amsart} \usepackage{amssymb,color,tikz}
\usetikzlibrary{matrix} \nonstopmode \textwidth=16.00cm

\usepackage{graphicx,color}

\usepackage{xcolor}

\numberwithin{equation}{section} 

\newcommand {\PP}{\mathbb{P}}
 \def\cocoa{{\hbox{\rm C\kern-.13em
      o\kern-.07em C\kern-.13em o\kern-.15em A}}}
\newtheorem{theorem}{Theorem}[section]
\newtheorem{lemma}[theorem]{Lemma}
\newtheorem{proposition}[theorem]{Proposition}
\newtheorem{corollary}[theorem]{Corollary}
 \theoremstyle{definition}
\newtheorem{definition}[theorem]{Definition} \theoremstyle{remark}
\newtheorem{remark}[theorem]{Remark}
\newtheorem{example}[theorem]{Example}

\DeclareMathOperator{\rank}{rank}
\DeclareMathOperator{\Ann}{Ann}

\DeclareMathOperator{\hess}{hess}
\DeclareMathOperator{\Hess}{Hess}

\usepackage{hyperref}
\definecolor{MyDarkGreen}{cmyk}{0.7,0,1,0}

\newcommand\Tstrut{\rule{0pt}{3ex}}
\newcommand\Bstrut{\rule[-1.6ex]{0pt}{0pt}}

\begin{document}

\title[Perazzo 3-folds and the weak Lefschetz property]
{Perazzo 3-folds and the weak Lefschetz property}

 \author[L. Fiorindo]{Luca Fiorindo} 
 \address{Dipartimento di Matematica e  Geoscienze, Universit\`a di
Trieste, Via Valerio 12/1, 34127 Trieste, Italy}
  \email{luca.fiorindo@dima.unige.it, ORCID 0000-0002-6435-0128}
  
 \author[E. Mezzetti]{Emilia Mezzetti} 
 \address{Dipartimento di Matematica e  Geoscienze, Universit\`a di
Trieste, Via Valerio 12/1, 34127 Trieste, Italy}
  \email{mezzette@units.it}
  \author[R.\ M.\ Mir\'o-Roig]{Rosa M.\ Mir\'o-Roig} 
  \address{Facultat de
  Matem\`atiques i Inform\`atica, Universitat de Barcelona, Gran Via des les
  Corts Catalanes 585, 08007 Barcelona, Spain} \email{miro@ub.edu,  ORCID 0000-0003-1375-6547}

\thanks{The second  author has been  partially supported by the grant PRIN $2017$SSNZAW$\_005$ \lq\lq Moduli Theory and Birational Classification'' and is member of INdAM - GNSAGA} 
\thanks{The third author has been partially supported by the grant PID2019-104844GB-I00}

\begin{abstract} We deal with Perazzo 3-folds in $\PP^4$, i.e. hypersurfaces $X=V(f)\subset \PP^4$ of degree $d$ defined by a homogeneous polynomial $f(x_0,x_1,x_2,u,v)=p_0(u,v)x_0+p_1(u,v)x_1+p_2(u,v)x_2+g(u,v)$, where  $p_0,p_1,p_2$ are algebraically dependent but linearly independent forms of degree $d-1$ in $u,v$, and $g$ is a form in $u,v$ of degree $d$. Perazzo 3-folds have vanishing hessian and, hence, 
the associated graded Artinian Gorenstein algebra $A_f$ fails the strong Lefschetz Property. In this paper, we determine the maximum and minimum Hilbert function of $A_f$ and we prove that if $A_f$ has maximal Hilbert function it fails the weak Lefschetz Property while  it satisfies the weak Lefschetz Property when it  has minimum Hilbert function. 
In addition, we classify all Perazzo 3-folds in $\PP^4$ such that $A_f$ has minimum Hilbert function.
\end{abstract}

\maketitle

\section{Introduction}
Having vanishing hessian is an elementary property of hypersurfaces that are cones. Conversely in every  projective space $\PP^N$ with $N\geq 4$ there exist classes of examples of hypersurfaces with vanishing hessian that are not cones. This was first proved by P. Gordan and M. Noether  (\cite{GN}), who disproved a claim by O. Hesse: any hypersurface $X\subset \PP^N$ with vanishing hessian is a cone  (\cite{He1}, \cite{He2}). 
Let $X$ be the hypersurface in $\PP^N$ defined by a polynomial $f(x_0,\ldots,x_N)$. Gordan and Noether realized that $X$ being a cone is equivalent to the condition that the partial derivatives of $f$ are $K$-linearly dependent, while $X$ has vanishing hessian if and only if they are $K$-algebraically dependent. 
They also gave a complete description of the hypersurfaces in $\PP^4$, not cones, with vanishing hessian.
Subsequent contributions were given by several authors. We refer to \cite{R} for an exhaustive bibliography.

J. Watanabe in \cite{w1}, and in \cite{MW} in collaboration with T.  Maeno, established the following strict connection between the Lefschetz properties of Artinian Gorenstein algebras and hypersurfaces with vanishing hessian. We recall that an Artinian $K$-algebra $A$ has the strong Lefschetz Property (respectively, the weak Lefschetz Property) if for a general linear form $L\in [A]_1$, the morphism $\times L^k:[A]_t\longrightarrow [A]_{t+k}$ has maximal rank for all integers $t\ge 0$ and $k\ge 1$ (respectively, the morphism $\times L:[A]_t\longrightarrow [A]_{t+1}$ has maximal rank for all integers $t\ge 0$).

Given a homogeneous polynomial $f$ in $N+1$ variables, we denote by $A_f$ the quotient of the differential operators' ring by the annihilator of $f$; it is a standard Artinian Gorenstein algebra whose socle degree $d$ coincides with the degree of $f$. 
In addition to the classical hessian, one defines the hessians of $f$ of order $t$, for $0\leq t\leq d$. Then, $A_f$ fails the strong Lefschetz Property if and only if the hessian of $f$ of order $t$ vanishes for some $t$ with $1\leq t\leq \lfloor{\frac{d}{2}}\rfloor$. 
In particular, the hypersurfaces with vanishing hessian all fail the strong Lefschetz Property. A natural question is then if they have or fail the weak Lefschetz Property. This question was considered by R. Gondim in \cite{G}, and he found examples of both types.

In this article, we consider the case of $\PP^4$, where the classification of hypersurfaces with vanishing hessian not cones is complete.
Following the terminology introduced by Gondim in \cite{G}, a hypersurface  in $\PP^4$ of degree $d\geq 3$ is a Perazzo hypersurface if, using homogeneous coordinates $x_0, x_1, x_2, u,v$, it has equation of the form $f=p_0x_0+p_1x_1+p_2x_2+g$, where  $p_0,p_1,p_2$ are  algebraically dependent but linearly independent forms of degree $d-1$ in $u,v$, and $g$ is a form in $u,v$ of degree $d$. Perazzo hypersurfaces have vanishing hessian. On the other hand, according to  \cite{Wa}, \cite{Wa1} and \cite{WB}, any hypersurface of degree $d$, with $3\leq d\leq 5$  of $\PP^4$ not cone with vanishing hessian is a Perazzo hypersurface.
In general, as proved by Gordan and Noether in \cite{GN}, all forms with vanishing hessian, not cones, are elements of  $K[u,v][\Delta]$ where $\Delta$ is a Perazzo polynomial of the form $p_0x_0+p_1x_1+p_2x_2$ (see \cite[Theorem 7.3]{WB}).
If $d=3$ for such an $f$ clearly $A_f$ fails the weak Lefschetz Property. For $d=4$ Gondim proved that the Artinian Gorenstein algebra of every Perazzo $3$-fold has the weak Lefschetz Property. 

Here we study the weak Lefschetz Property for Artinian Gorenstein algebras associated to Perazzo $3$-folds of any degree $d\geq 3$. 
First we consider the possible Hilbert functions $HF_{A_f}$ of $A_f$. In Propositions \ref{upper} and \ref{lower} we prove that they have a maximum and a minimum, coinciding only if $d=3,4$. 
We then study the weak Lefschetz Property for the algebras $A_f$ whose Hilbert function attains the upper or the lower bound. 
Our main results, contained in Theorems \ref{main1} and \ref{main2}, say that $A_f$ has the weak Lefschetz Property if $HF_{A_f}$ is minimum but it fails the weak Lefschetz Property if $HF_f$ is maximum. 
For $d=5$ this exhausts all possibilities; for $d\geq 6$  we give examples proving that for intermediate values of the Hilbert function both possibilities occur (see Example \ref{intermediate_cases}).  For further results on this topic see \cite{A}.

We then focus our attention on  Artinian Gorenstein algebras $A_f$ having minimum Hilbert function; using the theory of Waring rank for forms in $2$ variables, we are able to obtain in Theorem \ref{degree_d} a complete list  of these Perazzo $3$-folds. 
The classification is in terms of the position of the linear space $\pi$ generated by $p_0, p_1, p_2$ in $\PP(K[u,v]_{d-1})$, with respect to the secant varieties of the rational normal curve $C_{d-1}$. 
It results that, to ensure that $A_f$ has minimal Hilbert function,  $\pi$ has to meet $C_{d-1}$, and there are three possibilities: either $\pi$ is an osculating plane to $C_{d-1}$, or it is tangent to $\pi$ and meets the curve again in a second point, or the intersection $\pi\cap C_{d-1}$ consists of three distinct points. We conclude with a geometrical study of the polar and Gauss maps associated to these $3$-folds.
\par
\medskip

Next we outline the structure of this article. In Section \ref{prelim}, we recall the notions of strong and weak Lefschetz Property of an Artinian Gorenstein algebra, and of higher order hessians of a form. Then we state the theorem of J. Watanabe establishing a link between the failure of the strong Lefschetz Property for Artinian Gorenstein algebras and the vanishing of some hessian (Theorem \ref{watanabe}). We give some examples illustrating these notions.
In Section \ref{hvector}, we define Perazzo hypersurfaces and we study the $h$-vectors $h=(h_0,\ h_1, \ \ldots , h_{d-1}, \ h_d)$  of the associated Artinian Gorenstein algebras in the case of $\PP^4$. We relate them to the ranks of some block matrices composed of catalecticant matrices. In Propositions \ref{upper} and \ref{lower} we find the minimum and the maximum $h$-vector of these algebras for any degree $d\geq 4$. 
In Section \ref{wlp}, we study if the weak Lefschetz Property holds for Artinian Gorenstein algebras associated to Perazzo hypersurfaces.  We prove that, for $d\geq 5$, the algebras $A_f$ whose $h$-vector is maximum always fail the weak Lefschetz Property (Theorem \ref{main1}), while  the algebras whose $h$-vector is minimum always have it (Theorem \ref{main2}).
In Section \ref{classification}, we give a full classification of the Perazzo $3$-folds of degree $d\geq 5$ whose associated Artinian Gorenstein algebra has minimum $h$-vector. This is done using the  stratification of $\PP^{d-1}$ via the Waring rank of forms of degree $d-1$ in two variables. Finally we study the geometry of these $3$-folds, in terms of the image and fibres of their polar and Gauss map.

\vskip 4mm

\noindent {\bf Acknowledgement}. Most of this work was done while the third author was a guest of the University of
Trieste, and she would like to thank the people of  the Dipartimento di Matematica e  Geoscienze for their warm hospitality.

The authors wish to thank the anonymous referee for helpful comments.

\section{Lefschetz properties and Artinian Gorenstein algebras}\label{prelim}

In this section we fix notation, we recall the definition of weak/strong Lefschetz Property and we briefly discuss
general facts on Artinian Gorenstein algebras needed in next sections.

Throughout this work  $K$ will be an algebraically closed field of characteristic zero.
 Given a standard graded Artinian $K$-algebra $A=R/I$ where $R=K[x_0,x_1,\dots,x_n]$ and $I$ is a homogeneous ideal of $R$,
we denote by $HF_A:\mathbb{Z} \longrightarrow \mathbb{Z}$ with $HF_A(j)=\dim _K[A]_j$
its Hilbert function. Since $A$ is Artinian, its Hilbert function is
captured in its {\em $h$-vector} $h=(h_0,h_1,\dots ,h_e)$ where $h_i=HF_A(i)>0$ and $e$ is the last index with this property. The integer $e$ is called the {\em socle degree of} $A$.
\subsection{Lefschetz properties}

\begin{definition}
Let $A=R/I$ be a graded Artinian $K$-algebra. We say that $A$ has the {\em weak Lefschetz Property} (WLP, for short)
if there is a linear form $L \in [A]_1$ such that, for all
integers $i\ge0$, the multiplication map
\[
\times L: [A]_{i}  \longrightarrow  [A]_{i+1}
\]
has maximal rank, i.e.\ it is injective or surjective.
(We will often abuse notation and
say that the ideal $I$ has the WLP.) In this case, the linear form $L$ is called a Lefschetz
element of $A$. If for the general form $L \in [A]_1$ and for an integer $j$ the
map $\times L:[A]_{j-1}  \longrightarrow  [A]_{j}$ does not have maximal rank, we will say that the ideal $I$ fails the WLP in
degree $j$.

$A$ has the {\em strong Lefschetz Property} (SLP, for short) if there is a linear form $L \in [A]_1$ such that, for all
integers $i\ge0$ and $k\ge 1$, the multiplication map
\[
\times L^k: [A]_{i}  \longrightarrow  [A]_{i+k}
\]
has maximal rank.  Such an element $L$ is called a strong Lefschetz element for $A$.

$A$  has the {\em strong Lefschetz Property in the narrow sense} if there
exists an element $L \in [A]_1$  such that the multiplication map
\[
\times L^{e-2i}: [A]_{i}  \longrightarrow  [A]_{c-i}
\]
is bijective for $i=0,\cdots , [e/2]$ being $e$ the socle degree of $A$. 
\end{definition}

At first glance the problem of determining whether an Artinian standard graded $K$-algebra $A$ has the WLP seems a 
simple problem of linear algebra, but instead it has proven to be extremely elusive. Part
of the great interest in the WLP stems from the ubiquity of its presence and there  are a long series of papers determining classes of Artinian
algebras holding/failing the WLP but much more work remains to be done (see, for instance, \cite{CMM-R} and \cite{MMO}). The first result in this direction is due to Stanley \cite{s} and Watanabe \cite{w} and it asserts that the
WLP holds for an Artinian complete intersection ideal generated by powers of linear forms. 

\begin{example}
(1)
The ideal $I = (x_1^3,x_1^3,x_2^3,x_1x_2x_3)\subset  K[x_1,x_2,x_3]$ fails to have the WLP,
because for any linear form $L = ax_1 + bx_2 + cx_3$ the multiplication map
$$\times L: [k[x_1,x_2,x_3]/I]_2\cong K^6\longrightarrow  [k[x_1,x_2,x_3]/I]_3\cong K^6$$
is neither injective nor surjective.  More details on this example can be found in \cite[Example 3.1]{BK}.

(2) The ideal $I=(x_1^3,x_2^3,x_3^3,x_1^2x_2 )\subset K[x_1,x_2,x_3]$ has the WLP.
Since the $h$-vector of $R/I$ is (1,3,6,6,4,1), we only need to
check that  the map $\times L : [R/I]_{i} \longrightarrow [R/I]_{i+1}$
induced by $L=x_1+x_2+x_3$ is surjective for $i=2,3,4$. This is equivalent
to check that $[R/(I,L)]_{i}=0$ for $i=3,4,5$. Obviously, it is enough
to check the case $i=3$. We have

\[
\begin{array}{lll}
[R/(I,L)]_3 & \cong  [K[x_1,x_2,x_3]/(x_1^3,x_2^3,x_3^3,x_1^2x_2,x_1+x_2+x_3)]_3 \\ & \cong
[K[x_1,x_2]/(x_1^3,x_2^3,x_1^3+3x_1^2x_2+3x_1x_2^2+x_2^3, x_1^2x_2)]_3 \\ &
\cong  [k[x_1,x_2]/(x_1^3,x_2^3,x_1^2x_2,x_1x_2^2)]_3=0 \end{array}\]
which proves what we want.
\end{example}

It is worthwhile to point out that the weak Lefschetz
Property implies the unimodality of the Hilbert function.
If a graded Artinian $K$-algebra $A$ has the SLP in the
narrow sense, then the Hilbert function of A is unimodal and symmetric. Finally, 
if  a graded Artinian $K$-algebra $A$ has a symmetric Hilbert function, the
notion of the SLP on $A$ coincides with the one in the
narrow sense. In this work, we will deal with Artinian Gorenstein algebras $A$. It is well known that $A$ has symmetric Hilbert function. So, in the subsequent sections, the strong Lefschetz Property will be used in the narrow
sense.

\subsection{Artinian Gorenstein ideals}

In this subsection, we will characterize the Lefschetz elements for graded Artinian
Gorenstein algebras $A$. Given $R=K[x_0,\cdots ,x_n]$, we denote by $S=K[y_0,\cdots ,y_n]$ the ring of differential operators on $R$, i.e., $y_i=\frac{\partial}{\partial x_i}$. For any homogeneous polynomial $f\in R_d$, we define
$$\Ann_S(f):=\{ p\in S \mid p(f)=0\}\subset S.
$$
It is well known that $A=S/\Ann_S(f)$ is a standard graded Artinian Gorenstein $K$-algebra. Conversely, the theory of inverse systems developed by Macaulay  gives the following characterization of standard graded Artinian Gorenstein $K$-algebras.

\begin{proposition} Set $R=K[x_0,\cdots ,x_n]$ and let $S=K[y_0,\cdots ,y_n]$ be the ring of differential operators on $R$. Let $A=S/I$ be a standard  Artinian graded $K$-algebra. Then, $A$ is Gorenstein if and only if there is $f\in R$ such that $A\cong S/Ann_S(f)$. Moreover, isomorphic Gorenstein algebras are defined by forms equal up to a linear change of variables in $R$.
\end{proposition}

Under the hypothesis of the above proposition we have that the degree of $f$ coincides with the socle degree of $A$. 

\vskip 4mm
\begin{definition}
Let $f\in K[x_0, \cdots  , x_n]$ be a homogeneous polynomial and let $ A = S/\Ann_S(f)$ be the
associated Artinian Gorenstein algebra. Fix $\mathcal{B} = \{w_j \mid 1\le j \le h_t:=\dim A_t \} \subset A_t$ be an ordered
$K$-basis. The $t$-th (relative) {\em Hessian matrix} of $f$ with respect to $\mathcal{B}$ is defined as the $h_t \times h_t$ matrix:
$$
\Hess_f^t=(w_iw_j(f))_{i,j}.$$
The $t$-th {\em Hessian of } $f$ {\em with respect to} $\mathcal{B}$ is 
$$ \hess _f^t=\det (\Hess _f^t).
$$
\end{definition}

The 0-th Hessian is just the polynomial $f$ and, in the case $\dim A_1=n+1$, the 1st Hessian, with respect to the standard basis, is the classical Hessian. It is worthwhile to point out that the definition of Hessians and Hessian matrices of order $t$ depends on the choice of a basis of $A_t$ but  the vanishing of the $t$-th Hessian
is independent of this choice.

We end this preliminary section with a result due to Watanabe which establishes a useful link between the failure of Lefschetz properties and the vanishing of higher order Hessians.

\begin{theorem} \label{watanabe}
Let $f\in K[x_0, \cdots  , x_n]$ be a homogeneous polynomial of degree $d$ and let $ A = S/\Ann_S(f)$ be the
associated Artinian Gorenstein algebra. $L=a_0y_0+\cdots +a_ny_n\in A_1$ is a strong Lefschetz element of $A$ if and only if $\hess _f^t(a_0,\cdots, a_n)\ne 0$ for $t=1,\cdots,[d/2]$. More precisely, up to a multiplicative constant, $\hess _f^t(a_0,\cdots, a_n)$ is the determinant of the dual of the multiplication map 
$\times L^{d-2t}: [A]_{t}  \longrightarrow  [A]_{d-t}.$
\end{theorem}
\begin{proof}
See \cite[Theorem 4]{w1} and \cite[Theorem 3.1]{MW}.
\end{proof}

\begin{example} To illustrate Watanabe's theorem, we consider Ikeda's example of an Artinian Gorenstein algebra of codimension 4 failing WLP (see \cite[Example 4.4]{I}). We take $$f=x_0x_2^3x_3+x_1x_2x_3^3+x_0^3x_1^2\in K[x_0,x_1,x_2,x_3].$$  Let $S=K[y_0,y_1 ,y_2,y_3]$ be the ring of differential operators on $R$. We compute $\Ann _S(f)$ and we get:
$$\Ann_S(f)=\langle  y_0y_3^2, y_1^2y_3, y_0y_1y_3, y_0^2y_3, y_1y_2^2, y_0y_2^2-y_1y_3^2, y_1^2y_2 ,y_0y_1y_2, y_0^2y_2, y_1^3, y_3^4 , y_2^2y_3^2,$$ $$ y_2^4, y_0^2y_1^2-2y_2^3y_3,     y_0^3y_1-2y_2y_3^3, y_0^4 \rangle .$$ The $h$-vector of $A=S/\Ann_S(f)$ is: $(1, \ 4, \ 10, \ 10, \ 4, \ 1)$. We will apply the above criterion to check that $A$ fails the WLP in degree 3. To this end, we consider a $K$-basis of $[A]_2$: $$\{y_0^2, y_1^2, y_2^2, y_3^2, y_0y_1, y_0y_2, y_0y_3, y_1y_2, y_1y_3, y_2y_3 \}.$$ We get

$$
\Hess_f^2=\begin{pmatrix}
 0 & 12x_0 & 0 & 0 & 6x_1 & 0 & 0 & 0 &  0 & 0  \\
 12x_0 & 0 & 0  & 0 & 0 & 0 & 0 & 0 & 0 &  0  \\
  0 & 0 & 0 & 0 & 0 & 6x_3 & 6x_2 & 0 &  0 & 6x_0  \\
   0 & 0 & 0 & 0 & 0 & 0 & 0 &  6x_3 & 6x_2 & 6x_1 \\
   6x_1 & 0 & 0 & 0 & 6x_0 & 0 & 0 & 0 &  0 & 0  \\
   0 & 0 & 6x_3 & 0 & 0 & 0 & 0 &  0 &  0 & 6x_2  \\
   0 & 0 & 6x_2 & 0 & 0 & 0 & 0 & 0 &  0 & 0  \\
   0 & 0 & 0 & 6x_3 & 0 & 0 & 0 & 0 &  0 & 0  \\
    0 & 0 & 0 & 6x_2 & 0 & 0 & 0 & 0 &  0 & 6x_3  \\
     0 & 0 & 6x_0 & 6x_1 & 0 & 6x_2 & 0 & 0 &  6x_3 & 0  
\end{pmatrix}
.$$
For any $(a_0,a_1,a_2,a_3)\in K^4$, we have  $\hess _f^2(a_0,a_1,a_2,a_3)=0$. So, for any $L\in [A]_1$, the multiplication map $
\times L: [A]_{2}  \longrightarrow  [A]_{3}
$ has zero determinant. Therefore, it is not bijective and we conclude that $A$ fails the WLP. 
\end{example} 


\section{Perazzo 3-folds and the h-vector of the associated Gorenstein algebra}\label{hvector}

The goal of this section is to get upper and lower bounds for the $h$-vector of a standard graded Artinian Gorenstein algebra associated to a Perazzo 3-fold $X$ in $\PP^4$. So, let us start recalling its definition.

\begin{definition}\label{perazzo} Fix $N\ge 4$. A {\em Perazzo} hypersurface  $X\subset \PP^N$ of degree $d$ is a hypersurface defined by a form $f\in K[x_0,\cdots ,x_n,u_1\cdots ,u_m]$ of the following type:
$$ f=x_0p_0+x_1p_1+\cdots +x_np_n+g
$$
where $n+m=N$, $n,m\ge 2$, $p_i\in K[u_1,\cdots ,u_m]_{d-1}$ are algebraically dependent but linearly independent and $g\in K[u_1,\cdots ,u_m]_{d}$.
\end{definition}

It is worthwhile to point out that usually Perazzo hypersurfaces are assumed to be reduced and irreducible (see, for instance, \cite[Definition 3.12]{G}). We will insert these hypotheses if it is required.

\begin{example}
As a first example of Perazzo hypersurface we have the cubic 3-fold in $\PP^4$ of equation:
$$ f(x_0,x_1,x_2,u,v)=x_0u^2+x_1uv+x_2v^2.
$$
It is a cubic hypersurface with vanishing hessian but not a cone. So, it provides the first counterexample to Hesse's claim: any hypersuface $X\subset \PP^N$ with vanishing hessian is a cone (\cite{He1} and \cite{He2}).
\end{example}

\vskip 2mm
Hesse's claim, which is true for quadrics, was studied by Gordan and Noether in \cite{GN} for hypersurfaces of degree $d\geq 3$. They proved it is true for $N\leq 3$ but it is false for any $N\geq 4$. More precisely they gave a complete classification of the hypersurfaces with vanishing hessian for $N=4$ and a series of examples of hypersurfaces with vanishing hessian not cones for any $N\geq 5$. Subsequently Perazzo in \cite{P} described all cubic hypersurfaces with vanishing hessian for $N=4, 5, 6$.
The results of Gordan-Noether and of Perazzo have been recently considered and rewritten in modern language by many authors \cite{B}, \cite{CRS}, \cite{F}, \cite{L}, \cite{GR}, \cite{Wa} and \cite{WB}. 

\begin{remark}\label{linking_remark}
   We recall that the hypersurface defined by a polynomial $f$ has vanishing hessian if and only if the partial derivatives of $f$ are algebraically dependent, and it is a cone if and only if they are linearly dependent. It follows that the Perazzo hypersurfaces, introduced in Definition \ref{perazzo}, have all vanishing first hessian and in general are not cones.
\end{remark}

In $\PP^4$ the Gordan-Noether classification states that, for degree $d\leq 5$,  the hypersurfaces not cones with vanishing hessian are all Perazzo hypersurfaces, while for degree  $d> 5$, a form of degree $d$ with vanishing hessian, not cone, is an element of  $K[u,v][\Delta]$ where $\Delta$ is a Perazzo polynomial of the form $p_0x_0+p_1x_1+p_2x_2$ (see \cite{GN} and \cite[Theorem 7.3]{WB}).

In \cite{MW} Maeno and Watanabe found a connection between the vanishing of higher order hessians and Lefschetz properties, in particular with the SLP; then Gondim in \cite{G} studied the WLP for some hypersurfaces  with vanishing hessian.

In this note, we will concentrate our attention on Perazzo 3-folds $X$ in $\PP^4$ and our first goal will be to determine the  maximum and minimum $h$-vector for the Gorenstein Artinian algebras associated to  them. We will use the following notations: $R=K[x_0,x_1,x_2,u,v]$ is the polynomial ring in $5$ variables, $S=K[y_0,y_1,y_2,U,V]$ is the ring of differential operators on $R$, and a Perazzo $3$-fold $X\subset \PP^4$ of degree $d$ is defined by a form  
\begin{equation}\label{perazzo form}
f=x_0p_0(u,v)+x_1p_1(u,v)+x_2p_2(u,v)+g(u,v)\in R_d.\end{equation} 

If $d=3$, the corresponding algebras have all the same $h$-vector, and precisely $(1,5,5,1)$. In fact, by Remark \ref{linking_remark}, $X$ not being a cone implies $h_1=h_2=5$.  So, from now on, we will assume that $d\ge 4$ and we write
  \begin{equation}\label{coefficients perazzo}
      \begin{array}{rcl} p_0(u,v) & =& \sum _{i=0}^{d-1} {\binom{d-1} {i}}a_iu^{d-1-i}v^{i}, \\
p_1(u,v) & = & \sum _{i=0}^{d-1}{\binom{d-1} i}b_iu^{d-1-i}v^{i}, \\
p_2(u,v) & = & \sum _{i=0}^{d-1}{\binom{d-1} i}c_iu^{d-1-i}v^{i}, \text{ and } \\ g(u,v) & = & \sum _{i=0}^{d}{\binom{d} i}g_iu^{d-i}v^{i}. \end{array} \end{equation}
For any  $2\le k\le \lfloor \frac{d+1}{2} \rfloor$,  we define the matrices:

$$
\mathcal{A}_k:=\begin{pmatrix}
 a_0 & a_1 & \cdots & a_{k-1} \\
 a_1 & a_2 & \cdots & a_{k} \\
 \vdots & \vdots & & \vdots \\
  a_{d-k} & a_{d-k+1} & \cdots & a_{d-1} 
\end{pmatrix},
\quad 
\mathcal{B}_k:=\begin{pmatrix}
 b_0 & b_1 & \cdots & b_{k-1} \\
 b_1 & b_2 & \cdots & b_{k} \\
 \vdots & \vdots & & \vdots \\
  b_{d-k} & b_{d-k+1} & \cdots & b_{d-1} 
\end{pmatrix},  
$$

$$
\mathcal{C}_k:=\begin{pmatrix}
 c_0 & c_1 & \cdots & c_{k-1} \\
 c_1 & c_2 & \cdots & c_{k} \\
 \vdots & \vdots & & \vdots \\
  c_{d-k} & c_{d-k+1} & \cdots & c_{d-1} 
\end{pmatrix}, \text{ and } \quad \mathcal{G}_k:=
\begin{pmatrix}
    g_0 & g_1 & \cdots & g_k\\
    g_1 & g_2 & \cdots & g_{k+1}\\
    \vdots & \vdots & & \vdots \\
    g_{d-k} & g_{d-k+1} & \cdots & g_d
\end{pmatrix}.
$$
The matrices $\mathcal{A}_k$, $\mathcal{B}_k$, $\mathcal{C}_k$ and $\mathcal{G}_k$ are the building blocks of the  matrices $M_k$, $N_k$ and $N'_k$ that will play an important role in the proof of our main results. They are defined as follows:

$$
M_k:= \begin{pmatrix} \mathcal{A}_k | \mathcal{B}_k | \mathcal{C}_k
\end{pmatrix}, \quad \
 \quad N_k:=\begin{pmatrix}
    \mathcal{A}_{k+1} \\ \hline \mathcal{B}_{k+1} \\ \hline \mathcal{C}_{k+1}
\end{pmatrix} \text{ and } \quad N'_k:= \begin{pmatrix}
    \mathcal{A}_{k+1} \\ \hline \mathcal{B}_{k+1} \\ \hline \mathcal{C}_{k+1} \\ \hline \mathcal{G}_k
\end{pmatrix}.
$$

\begin{remark}
(1) The matrices $N_k$ and $M_{k+1}$ contain the same 3 blocks of size $(d-k)\times (k+1).$

(2) Since $M_k=N_{d-k}^t$, we have $\rank M_k=\rank N_{d-k}$.

(3) As we will see in the proof of Propositions \ref{upper} and \ref{lower}, the $h$-vector of $S/\Ann_S(f)$ is minimal if and only if for all $k$, $2\le k \le   \lfloor \frac{d}{2} \rfloor$, $\rank M_k=\rank N'_k=3$.
\end{remark}

\begin{proposition}\label{hilbert function}
Let $f=x_0p_0(u,v)+x_1p_1(u,v)+x_2p_2(u,v)$ be a form of degree $d$ defining a Perazzo $3$-fold in $\PP^4$. Let $h=(h_0, h_1, \ldots, h_d)$ be its $h$-vector. Then
$h_0=h_d=1, h_1=h_{d-1}=5$ and, for $2\leq i\leq d-2$, $h_i=4i+1-m_i-n_i$, where $m_i=3i-\rank M_i$ and  $n_i=i+1-\rank N_i$. 
\end{proposition}
\begin{proof} Recall  that the $h$-vector of an Artinian Gorenstein algebra is symmetric and, hence, we only have to compute $h_i$ for $0\le i \le \lfloor \frac{d}{2} \rfloor $.
We have $$h_i=\dim A_i=\dim S_i-\dim \Ann_S(f)_i={\binom{4+i}{i}}-\dim \Ann_S(f)_i. $$ So, we have to compute $\dim \Ann_S(f)_i$ for any $i$,  $0\le i \le \lfloor \frac{d}{2} \rfloor $. Since $p_0(u,v)$, $p_1(u,v)$ and $p_2(u,v)$ are $K$-linearly independent, we have $\dim \Ann _S(f)_1=0$ and, hence, $h_1=5$.  

We observe that, for any $i\geq 2$, $\Ann_S(f)_i$ contains $(y_0,y_1,y_2)^{i-k}(U,V)^k$, for $0\leq k\leq i-2$. Therefore 
$$\dim A_i\leq {\binom{4+i}{i}}-\sum_{k=0}^{i-2}(k+1){\binom{i-k+2}{2}}=4i+1.$$
We have to compute the numbers
$$m_i=\dim(\Ann_S(f)_i\cap (y_0,y_1,y_2)(U,V)^{i-1}),$$
$$n_i=\dim(\Ann_S(f)_i\cap (U,V)^i),$$
and we will get 
\begin{equation}\label{hi}
    \dim A_i=4i+1-m_i-n_i.
    \end{equation}
This can be done because there are no linear dependence relations between the two parts, given the bi-homogeneous nature of $f$ with respect to the two groups of variables $x_0,x_1,x_2$ and $u,v$. 

To compute $m_i$ we consider a general polynomial of degree $i$ in $(y_0,y_1,y_2)(U,V)^{i-1}$: 
$$(\alpha_0U^{i-1}+\alpha_1U^{i-2}V+\cdots+\alpha_{i-1}V^{i-1})y_0+(\beta_0U^{i-1}+\cdots+\beta_{i-1}V^{i-1})y_1+(\gamma_0U^{i-1}+\cdots)y_2.$$
It belongs to $\Ann_S(f)_i$ if and only if
$$\alpha_0p_{0,u^{i-1}}+\alpha_1p_{0,u^{i-2}v}+\cdots+\alpha_{i-1}p_{0,v^{i-1}}+\beta_0p_{1,u^{i-1}}+\cdots+\gamma_0p_{2,u^{i-1}}+\cdots +\gamma _{i-1}p_{2,v^{i}}=0.
$$
The partial derivatives of $p_0, p_1, p_2$ appearing in the above expression have degree $d-i$; setting equal to zero the coefficients  of the $d-i+1$ monomials in $u,v$, we get a homogeneous linear system of $d-i+1$ equations in the $3i$ unknowns $\alpha_0,\ldots,\alpha_{i-1}, \beta_0, \ldots, \beta_{i-1}, \gamma_0, \ldots, \gamma_{i-1}.$ The matrix of the coefficients is $M_i$, therefore $m_i=3i-\rank M_i$, and we are done.

To compute $n_i$ we consider a general polynomial of degree $i$ in $U,V$:
$$\delta_0U^i+\delta_1U^{i-1}V+\cdots+\delta_iV^i$$
and we impose that it belongs to $\Ann_S(f)_i$. We get
$$(\delta_0p_{0,u^i}+\delta_1p_{0,u^{i-1}v}+\cdots+\delta_ip_{0,v^i})x_0+(\delta_0p_{1,u^i}+\cdots)x_1+(\delta_0p_{2,u^i}+\cdots)x_2=0. $$

Looking at the coefficients of $x_0, x_1, x_2$ and then the coefficients of the monomials in $u,v$ of degree $d-i-1$ and $d-i$, we get a homogeneous linear system of $3(d-i)+(d-i+1)$ equations in $i+1$ unknowns, whose matrix of the coefficients is $N_i$. We conclude that $n_i=i+1-\rank N_i.$ The proof is complete.
\end{proof}

\begin{remark}\label{hilbert function bound}
We observe that the expression for $h_i$ can also be written in the form $h_i=\rank M_i+\rank N_i$. In fact, if we write a unique linear system to compute the dimension of the space $\Ann_S(f)_i\cap [(y_0, y_1, y_2)(U,V)^{i-1}+(U,V)^i]$,  the matrix of this linear system results to be $\left( \begin{array}{c|c}
   0 & N_i \\
   \hline
   M_i & 0 \\
\end{array}\right).$
  
  In the general case, when $f$ is as in (\ref{perazzo form}) with $g\ne 0$, equality (\ref{hi}) is not necessarily true, but only the inequality $\dim A_i\geq 4i+1-m_i-n_i$ holds true. An explicit example is provided by the form $f=x_0u^9+x_1u^8v+x_2v^9+u^5v^5$.
  
  On the other hand, in this more general situation the matrix associated to the linear system to be considered to compute $h_i$ is 
   $
   \left( \begin{array}{c|c}
   0 & N_i \\
   \hline
   M_i & \mathcal{G}_i \\
\end{array}\right)$.
   This implies, for every index $i$, the series of inequalities
   $$\rank M_i+\rank N_i\le h_i\le \rank M_i+\rank N'_i.$$
   
   Clearly, every time $\rank N_i=\rank N'_i$, we obtain a relation as in Proposition \ref{hilbert function}. This is obviously the case when $g=0$. It is also the case if one of the  polynomials $p_0$, $p_1$, $p_2$ is general enough. Indeed, we observe that $N_i$ has maximal rank if and only if its columns are linearly independent; so if the rank of one of the matrices $\mathcal{A}_{i+1}$, $\mathcal{B}_{i+1}$, $\mathcal{C}_{i+1}$ is computed by the number of columns then $N_i$ has maximal rank. This happens if one of the  polynomials $p_0$, $p_1$, $p_2$ is general enough in view of \cite[Proposition 3.4]{Iar}.
\end{remark}

\begin{proposition}\label{upper}  Let $d\ge 4$. The maximum $h$-vector of  the Artinian Gorenstein algebras $S/Ann_S(f)$ associated to the Perazzo $3$-folds of degree $d$ in $\PP^4$ is:
    \begin{itemize}
    \item[(1)]  If $d=4t-1$ 
    then 
        $$h_i= \begin{cases} 4i+1 \text{ for } 0\le i\le t \\
        4t+1 \text{ for } t+1\le i \le 2t-1 \\
        \text{ symmetry};
         \end{cases}$$
    \item[(2)] If  $d=4t$  
    then 
        $$h_i= \begin{cases} 4i+1 \text{ for } 0\le i \le t \\
        4t+2 \text{ for } t+1\le i \le 2t \\
        \text{ symmetry};
         \end{cases}$$
    \item[(3)] If $d=4t+1$ 
    then 
    $$h_i= \begin{cases} 4i+1 \text{ for } 0\le i\le t \\
    4t+3 \text{ for } t+1\le i \le 2t \\
    \text{ symmetry};
     \end{cases}
    $$
    \item[(4)] If $d=4t+2$   
    then 
        $$h_i= \begin{cases} 4i+1 \text{ for } 0\le i\le t \\
        4t+4  \text{ for } t+1\le i \le 2t+1 \\
        \text{ symmetry}.
         \end{cases}$$
    \end{itemize}
    \end{proposition}

\begin{proof} 
Let $f$ be a form of degree $d$ as in (\ref{perazzo form}) with $p_0, p_1, p_2, g$ as in (\ref{coefficients perazzo}). Being the $h$-vector symmetric, we only have to compute $h_i$ for $0\le i\le \frac{d}{2}$.

In view of Proposition \ref{hilbert function} and Remark \ref{hilbert function bound}, the maximal Hilbert function is obtained when $m_i, n_i$ are minimal for any $i$, i.e. when the ranks of the matrices $M_i$, $N'_i$ are as large as possible. 

Clearly $\rank M_i\leq \min\{3i, d-i+1\}.$ Therefore 
$$\rank M_i\leq \begin{cases} 3i \ \text{for} \  i\leq {\frac{d+1}{4}};\\
d-i+1 \ \text{for} \ i\geq {\frac{d+1}{4}}.
\end{cases}$$

Regarding $N'_i$, we observe that, in our situation,  $i+1\leq 3(d-i)+(d-i+1)$, so always $\rank N'_i\leq i+1$.

This gives upper bounds on $h_i$ depending on the class of congruence  of the degree $d$ modulo $4$, that are precisely those in the statement of this Proposition. 

To conclude the proof, we claim  that these bounds are achieved. To this end, we observe that, in view of the expressions (\ref{coefficients perazzo}), the columns of the matrices $\mathcal{A}_i, \mathcal{B}_i, \mathcal{C}_i, \mathcal{G}_i$ contain up to a constant the coefficients of the partial derivatives of order $i-1$ of $p_0, p_1, p_2, g$ respectively. But, if $p_0,p_1,p_2, g$ are general enough, then, by \cite[Proposition 3.4]{Iar}, for any $i$ their partial derivatives of order $i-1$  are as linearly independent as possible in $K[u,v]_{d-1-i}$. This means that the ranks of the matrices $\mathcal{A}_i, \mathcal{B}_i, \mathcal{C}_i, \mathcal{G}_i$ are as large as possible. This proves our claim.
\end{proof}

A class of explicit examples of polynomials such the bound in Proposition \ref{upper} is attained can be found in \cite[Example 3.20]{Fio}.

\vskip 2mm
To determine the minimum $h$-vector for the Gorenstein Artinian algebra associated to a Perazzo 3-fold $X$ in $\PP^4$ we need first to recall some results about the growth of the Hilbert function of standard graded $K$-algebras and to fix some additional notation.

Given  integers $n, d\ge 1$, we define the  {\em $d-$th binomial expansion of $n$} as

\vskip 2mm
$$
n=\binom{\eta_d}{d}+\binom{\eta_{d-1}}{d-1}+\cdots +\binom{\eta_e}{e}
$$

\noindent where $\eta_d>\eta_{d-1}>\cdots >\eta_e\ge e\ge 1$ 
are uniquely determined integers (see \cite[Lemma 4.2.6]{BH93}).
We write
\vskip 2mm
$$ 
n^{<d>}=\binom{\eta_d+1}{d+1}+\binom{\eta_{d-1}+1}{d}+\cdots +\binom{\eta_e+1}{e+1}, \text{ and }
$$
\vskip 2mm
$$
n_{<d>}=\binom{\eta_d-1}{d}+\binom{\eta_{d-1}-1}{d-1}+\cdots +\binom{\eta_e-1}{e}.
$$

\vskip 2mm
The numerical functions $H:\mathbb N \longrightarrow \mathbb N$ that are Hilbert functions of graded standard $K$-algebras were characterized by Macaulay, \cite{BH93}. Indeed, given a numerical function $H:\mathbb N \longrightarrow \mathbb N$ the following conditions are equivalent:
\begin{enumerate}
    \item[(i)] There exists a standard graded $K$-algebra $A$ with Hilbert function  $H$,
    \item[(ii)]  $H(0)=1$ and $H(t+1)\le H(t)^{<t>}$ for all $t \ge 1$.
\end{enumerate}

Notice that condition (ii) imposes strong restrictions on the Hilbert function of a standard graded $K$-algebra and, in particular, it bounds its growth. As an application of Macaulay's theorem, we have:

\begin{proposition}\label{lower} Let $d\geq 4$. Let $R=K[x_0,x_1,x_2,u,v]$ and $S=K[y_0,y_1,y_2,U,V]$ be the ring of differential operators on $R$. The minimum $h$-vector of  the
 Artinian Gorenstein algebras $A=S/Ann_S(f)$ associated to the Perazzo $3$-folds of degree $d$ in $\PP^4$ is:
$$ (1, \ 5, \ 6, \ 6, \ \cdots \ 6, \ 6, \ 5, \ 1).
$$
\end{proposition}
\begin{proof} The proof proceeds as follows: we first prove  that the cited $t$-uple is less than any possible $h$-vector associated to a Perazzo $3$-fold, with respect to the termwise order; then, we give examples of Perazzo forms that have this $h$-vector. Let  $$h_A=(h_0, \ h_1, \ h_2, \ h_3, \   \cdots \ h_{d-2}, \ h_{d-1}, \ h_d)$$ be the $h$-vector of $A$. First of all we observe that, arguing as in the proof of Proposition \ref{upper}, we get that $6\le h_2\le 9$ which, together with the fact that the $h$-vector of any standard graded Artinian Gorenstein algebra is symmetric, gives us that a lower bound for $h_A$
looks like $$(1, \ 5, \ 6, \ h_3, \   \cdots \ h_{d-2}, \ 6, \ 5, \ 1).$$
This concludes the first step for $d\le 5$. We will now prove that if $d\geq 6$, for any $i$, $3\le i \le d-3$, $h_i\ge 6$.

First we assume $d\ge 8$. If $h_j\le 5$ for some $5\le j \le d-3$, using Macaulay's inequality $h_{t+1}\le h_t^{<t>}$ for all $t\ge 1$, we get that $h_i\le 5$ for all $i\ge j$ contradicting the fact that $h_{d-2}=6$. Therefore, $h_j\ge 6$ for all $5\le j \le d-2$ and, by symmetry, we also have $h_3,h_4\ge 6.$ 

For $d=6,7$, we must show that $h_3\ge 6$. This last equality follows after a straightforward computation which shows that $(y_0,y_1,y_2)^3\oplus (y_0,y_1,y_2)^2(U,V)\subset \Ann _S(f)_3$, $\dim 
\{(\alpha _0U^2+\alpha _1UV+\alpha _2V^2)y_0+(\beta _0U^2+\beta _1UV+\beta _2V^2)y_1+(\gamma _0U^2+\gamma _1UV+\gamma _2V^2)y_2 \in \Ann _S(f)_3\}\le 6$ and  $\dim 
\{\delta _0U^3+\delta _1U^2V+\delta _2UV^2+\delta _3 V^3 \in \Ann _S(f)_3\}\le 1$. Therefore, $h_3=\dim S_3/\Ann_S(f)_3\ge {\binom{7}{4}}-29=6.$

Summarizing, we have got that for any $d\ge 4$ and for  $2\le i \le d-2$, it holds $h_i\ge 6$. To finish the proof  it suffices to give an example with   $h_i= 6$, for any $i$, $2\le i \le d-2$. We take the homogeneous polynomial of degree $d$:
$$f(x_0,x_1,x_2,u,v)=u^dx_0+u^{d-1}vx_1+v^dx_2.
$$
It is easy to check that it has the desired $h$-vector. 
\end{proof}

\begin{remark}\label{rank equal 3}
From Proposition \ref{hilbert function} it follows that $A=S/\Ann_S(f)$ has minimum $h$-vector $(1,5,6,6,\ldots,6, 6, 5, 1)$ if and only if $\rank M_i=\rank N'_i=3$ for any $i$ with $2\leq i\leq \lfloor\frac{d}{2} \rfloor.$ We note that none of these ranks can be strictly less than $3$ due to the assumption that $p_0, p_1, p_2$ are linearly independent.
\end{remark}

\begin{remark}
   From Propositions \ref{upper} and \ref{lower}, it follows that for $d=4$ the unique possible $h$-vector is $(1,5,6,5,1)$. Instead, for $d=5$, we can obtain only the maximal $h$-vector $(1,5,7,7,5,1)$, and the minimal $h$-vector $(1,5,6,6,5,1)$. For bigger values of $d$, also some intermediate cases are a priori possible.
\end{remark}

\section{Perazzo 3-folds and the WLP}\label{wlp}

From Theorem \ref{watanabe} and Remark \ref{linking_remark}, since the  Perazzo $3$-folds have vanishing first hessian, it follows that the associated algebras $A$ fail the strong Lefschetz Property. In particular the map
\[
\times L^{d-2}: [A]_{1}  \longrightarrow  [A]_{d-1}
\]
is not an isomorphism for every $L\in [A]_1$. 
The goal of this section is to analyze whether the Artinian  Gorenstein algebra $A$ associated to a Perazzo 3-fold $X\subset \PP^4$ has the WLP. If $d=3$, clearly $A$ fails also WLP. But Gondim has proved that, for any Perazzo $3$-fold of degree $4$, $A$  has the WLP (\cite{G}, Theorem 3.5).  More precisely, we will see that, in any degree $d\geq 5$, WLP holds when $A$ has minimum $h$-vector and fails when it has maximum $h$-vector.

\begin{theorem} \label{main1} Let $X\subset \PP^4$ be a Perazzo 3-fold of degree $d\ge 5$ and equation $$f=x_0p_0(u,v)+x_1p_1(u,v)+x_2p_2(u,v)+g(u,v)\in R_d=K[x_0,x_1,x_2,u,v]_d.$$ Let $S=K[y_0,y_1,y_2,U,V]$ be the ring of differential operators on $R$. If $A=S/\Ann_S(f)$ has maximum h-vector, then $A$ fails WLP.
\end{theorem}

\begin{proof} According to the parity of the socle degree of $A$, we distinguish two cases.

\noindent {\bf Case 1:} $d$ is odd. Write $d=2r+1$.To show that $A$ fails WLP, we will  prove that for any $L\in [A]_1$,  the multiplication map
\[
\times L: [A]_{r}  \longrightarrow  [A]_{r+1}
\]
is not bijective. By Theorem \ref{watanabe}, it is enough to see the vanishing of the $r$-th Hessian $\hess _f^r$ of $f=x_0p_0(u,v)+x_1p_1(u,v)+x_2p_2(u,v)+g(u,v)$ with respect to a suitable $K$-basis $\mathcal{B}$ of $[A]_r$.
First we can notice that a basis $\mathcal{B}$ made of classes with a monomial representative always exists. So, $\Hess _f^r$ is just a submatrix of dimension $h_r\times h_r$ of the following matrix:
$$\Big( \frac{\partial^{2r} f}{\partial u^\alpha\partial v^\beta\partial x_0^\gamma\partial x_1^\delta\partial x_2^\eta}\Big)_{\alpha+\beta+\gamma+\delta+\eta=2r}$$
where monomials are lexicographic ordered (for simplicity). Knowing that $f$ is linear in the variables $x_0,x_1,x_2$, the above matrix can be partially computed as:

\begin{equation*}
    \left(
    \begin{array}{cccc|cccc|ccc}
        \frac{\partial^{2r} f}{\partial u^{2r}} & \frac{\partial^{2r} f}{\partial u^{2r-1}\partial v} & \cdots  & \frac{\partial^{2r} f}{\partial u^{r}\partial v^{r}} & \frac{\partial^{2r-1} p_0}{\partial u^{2r-1}} &\frac{\partial^{2r-1} p_0}{\partial u^{2r-2}\partial v} & \cdots & \frac{\partial^{2r-1} p_2}{\partial u^{r}\partial v^{r-1}} & 0 & \cdots & 0 \Bstrut{}\\
        \frac{\partial^{2r} f}{\partial u^{2r-1}\partial v} & \frac{\partial^{2r} f}{\partial u^{2r-2}\partial v^2} & \cdots  & \frac{\partial^{2r} f}{\partial u^{r-1}\partial v^{r+1}} &\frac{\partial^{2r-1} p_0}{\partial u^{2r-2}\partial v} &\frac{\partial^{2r-1} p_0}{\partial u^{2r-3}\partial v^2} & \dots & \frac{\partial^{2r-1} p_3}{\partial u^{r-1}\partial v^{r}} & 0 & \cdots & 0 \Bstrut{}\\
        \vdots & \vdots & \ddots  & \vdots & \vdots & \vdots & \ddots  & \vdots & \vdots & \ddots & \vdots\Bstrut{}\\
        \frac{\partial^{2r} f}{\partial u^{r}\partial v^r} & \frac{\partial^{2r} f}{\partial u^{r-1}\partial v^{r+1}} & \cdots  & \frac{\partial^{2r} f}{\partial u^{r-1}\partial v^{2r}} & \frac{\partial^{2r-1} p_0}{\partial u^{r-1}\partial v^r} &\frac{\partial^{2r-1} p_0}{\partial u^{r-2}\partial v^{r+1}} & \dots & \frac{\partial^{2r-1} p_2}{\partial u^{2r-1}} & 0 & \cdots & 0\Bstrut{}\\
        \hline
        \frac{\partial^{2r-1} p_0}{\partial u^{r-1}\partial v^r} &\cdots & \cdots & \frac{\partial^{2r-1} p_0}{\partial u^{r-1}\partial v^r}  & 0 & \cdots & \cdots & 0 & 0 & \cdots & 0 \Tstrut{}\Bstrut{}\\
        \vdots & \vdots & \ddots  & \vdots & \vdots & \vdots & \ddots  & \vdots & \vdots & \ddots & \vdots \Bstrut{}\\
        \frac{\partial^{2r-1} p_2}{\partial u^r\partial v^{r-1}} &\cdots & \cdots & \frac{\partial^{2r-1} p_2}{\partial u^{2r-1}}  & 0 & \cdots & \cdots & 0 & 0 & \cdots & 0 \Bstrut{}\\
        \hline
        0 & \cdots & \cdots  & 0 & 0 & \cdots & \cdots  & 0 & 0 & \cdots & 0\Tstrut{}\\
        \vdots & \ddots & \ddots  & \vdots & \vdots & \ddots & \ddots  & \vdots & \vdots & \ddots & \vdots\\
        0 & \cdots & \cdots  & 0 & 0 & \cdots & \cdots  & 0 & 0 & \cdots & 0\\
    \end{array}
    \right)
\end{equation*}

\vskip 2mm
The three vertical (respectively, horizontal) blocks  are composed respectively by $r+1,\, 3r, \, \binom{r+4}{4}-(4r+1)$ columns (respectively, rows). Thus every possible choice of a $h_r\times h_r$ submatrix turns out to have at least an all zero sub-submatrix of size $(h_r-(r+1))\times (h_r-(r+1))$. We now use the hypothesis of $A$ to have maximum $h$-vector and Proposition \ref{upper} to obtain that $h_r=2r+3$. We have just proved that $\Hess _f^r$, matrix of dimension $(2r+3)\times (2r+3)$, has an all zero submatrix of dimension $(r+2)\times(r+2)$: this implies that $\hess _f^r$ identically vanishes.

\vskip 2mm

\noindent {\bf Case 2:} $d$ is even. Write $d=2r+2$. Note that, since the $h$-vector is maximum, then $h_r=h_{r+1}=h_{r+2}$. Using again the hessian criterion of Watanabe's (Theorem \ref{watanabe}), we will check that for any $L\in [A]_1$,  the multiplication map
\[
\times L^2: [A]_{r}  \longrightarrow  [A]_{r+2}
\]
is not bijective. This implies that for any $L\in [A]_1$,  the multiplication map
\[
\times L: [A]_{r}  \longrightarrow  [A]_{r+1}
\]
is not bijective and, hence, $A$ fails the WLP. 

Same adapted argument of the previous case can be used also here. In fact, the matrix to be considered is $\Hess _f^r$ which is now of size $(2r+4)\times(2r+4)$ which is even bigger than the previous case. Thus, as discussed above, its determinant is always zero.
\end{proof}

In contrast with the last result we have  that if an Artinian Gorenstein algebra $A$ associated to a Perazzo 3-fold has  minimum $h$-vector, then $A$ has the WLP. Our proof uses Green's theorem that we recall for sake of completeness.

\begin{theorem} \label{green} Let $A=R/I$ be an Artinian graded algebra and let $L\in A_1$ be a general linear form. Let $h_t$ be the entry of degree $t$ of the h-vector of $A$. Then the degree $t$ entry $h'_t$ of the h-vector of $R/(I,L)$ satisfies the inequality:
$$
h'_t\le (h_t)_{<t>} \text{ for all } t\ge 1.
$$
\end{theorem}
\begin{proof}
See \cite[Theorem 1]{Gr}.
\end{proof}

\begin{theorem} \label{main2} Let $X\subset \PP^4$ be a Perazzo 3-fold of degree  $d\ge 5$ and equation $$f=x_0p_0(u,v)+x_1p_1(u,v)+x_2p_2(u,v)+g(u,v)\in R=K[x_0,x_1,x_2,u,v]_d.$$ Let $S=K[y_0,y_1,y_2,U,V]$ be the ring of differential operators on $R$. If $A=S/\Ann_S(f)$ has minimum h-vector, then $A$ has WLP.
\end{theorem}

\begin{proof}
For $5\le d\le 7$ see next section where a full classification of Perazzo 3-folds with minimal $h$-vector is given. Assume $d\ge 8$. By the minimality assumption, $h_2=h_3=\cdots=h_{d-2}=6$. By \cite[Proposition 2.1]{MMN}, if for a general linear form $L\in [A]_1$, the multiplication map \[
\times L: [A]_{2}  \longrightarrow  [A]_{3}
\]
is bijective, then 
\[
\times L: [A]_{1}  \longrightarrow  [A]_{2}
\]
is injective, and for all $j\ge 2$,
\[
\times L: [A]_{j}  \longrightarrow  [A]_{j+1}
\]
is surjective, therefore $A$ has the WLP.  By the symmetry property of Artinian Gorenstein algebras, 
 \[
\times L: [A]_{2}  \longrightarrow  [A]_{3}
\]
is bijective if and only if
\[
\times L: [A]_{d-3}  \longrightarrow  [A]_{d-2}
\]
is bijective. So, let us prove the bijection of this last map. To this end, for a general linear form $L\in [A]_1$, we consider the exact sequence:
\[
 [A]_{d-3}  \longrightarrow  [A]_{d-2} \longrightarrow  [S/(\Ann _S(f),L)]_{d-2} \longrightarrow 0.
\] It follows that $\times L: [A]_{d-3}  \longrightarrow  [A]_{d-2}$ is bijective if and only if $[S/(\Ann _S(f),L)]_{d-2}=0$. Using the hypothesis $
d-2\ge 6$ (and, hence, $h_{d-2}\le d-2$) and Theorem \ref{green} we get
$$ \dim [S/(\Ann _S(f),L)]_{d-2} \le (h_{d-2})_{<d-2>}=0
$$
which proves what we want.
\end{proof}

\begin{remark}
   As a consequence of Theorem \ref{main2}, all forms of degree $d$ which define a Perazzo $3$-fold with minimum $h$-vector are examples of forms with zero first order hessian, and all hessians of order $t$ different from zero, for  $2\le t\le\lfloor\frac{d}{2}\rfloor$.
\end{remark}

For Gorenstein Artinian algebras associated to Perazzo 3-folds $X$ in $\PP^4$ and with intermediate $h$-vector both possibilities occur: there are examples failing WLP and examples satisfying WLP as next example shows.

\begin{example}\label{intermediate_cases}
1.-  Let $X\subset \PP^4$ be the  Perazzo 3-fold of equation
$$ f(x_0,x_1,x_2,u,v)=u^6x_0+(u^2v^4+u^4v^2)x_1+v^6x_2\in K[x_0,x_1,x_2,u,v]_7.
$$
 Let $S=K[y_0,y_1,y_2,U,V]$ be the ring of differential operators on $R$. We have
 $$\Ann _S(f)=\langle y_0^2,y_1^2,y_2^2,y_0y_1,y_0y_2,y_1y_2,y_0V, y_2U,  y_0U^2+15y_1U^2-15y_1V^2-y_2V^2,$$
 $$U^3V-UV^3, 15y_1U^4-y_2V^4, UV^5, V^7, U^7 \rangle .$$
Therefore, the Artinian Gorenstein algebra  $A=S/Ann_S(f)$ has $h$-vector: $(1, \ 5, \ 7, \ 8, \ 8, \ 7, \ 5, \ 1)$. Using Macaulay2 \cite{M2} we check that for a general linear form $L\in [A]_1$, the multiplication map  
 \[
\times L: [A]_{3}  \longrightarrow  [A]_{4}
\]
is bijective and, hence, $A$ satisfies the WLP. It does not have the SLP because for any linear form $L\in [A]_1$
 \[
\times L^3: [A]_{2}  \longrightarrow  [A]_{5}
\]
is not surjective.

2.- Let $X\subset \PP^4$ be the  Perazzo 3-fold of equation
$$ f(x_0,x_1,x_2,u,v)=u^6x_0+u^3v^3x_1+v^6x_2\in K[x_0,x_1,x_2,u,v]_7.
$$
 Let $S=K[y_0,y_1,y_2,U,V]$ be the ring of differential operators on $R$. We have
 $$\Ann _S(f)=\langle y_0^2,y_1^2,y_2^2,y_0y_1,y_0y_2,y_1y_2,
 y_0v, y_2u,  20y_1U^3-y_2V^3,$$
 $$y_0U^3-20y_1V^3, UV^4, U^4V, V^7, U^7 \rangle .$$
Therefore, the Artinian Gorenstein algebra  $A=S/Ann_S(f)$ has $h$-vector: $(1, \ 5, \ 7, \ 9, \ 9, \ 7, \ 5, \ 1)$. 
Computing the third hessian, since it results to be zero, we get
 that for any linear form $L\in [A]_1$, the multiplication map  
 \[
\times L: [A]_{3}  \longrightarrow  [A]_{4}
\]
is not bijective and, hence, $A$ fails the WLP.
\end{example}

\section{On the classification of certain Perazzo 3-folds of degree at least 5}\label{classification}

The goal of this section is to classify all Perazzo 3-folds $X$ in $\PP^4$ of degree $d\ge 5$ whose associated Artinian Gorenstein algebra $S/\Ann _S (f)$  has $h$-vector: $(1, \ 5, \ 6, \ 6, \cdots , \ 6, \ 6, \ 5, \ 1)$.
As a corollary we will also classify all Perazzo 3-folds $X$ in $\PP^4$ of degree 5 whose associated Artinian Gorenstein algebra $S/\Ann _S (f)$ has the WLP.

We start the section with some technical lemmas and remarks.
\begin{lemma}\label{generators} Let $f_1=p_0(u,v)x_0+p_1(u,v)x_1+p_2(u,v)x_2$ and $f_2=q_0(u,v)x_0+q_1(u,v)x_1+q_2(u,v)x_2$ be two Perazzo 3-folds of degree $d$ in $\PP^4$ such that $\langle p_0,p_1,p_2 \rangle = \langle q_0,q_1,q_2 \rangle \subset K[u,v]_{d-1} $.  Then, the $h$-vectors of $S/\Ann _S (f_1)$ and $S/\Ann _S (f_2)$ coincide.
\end{lemma}
\begin{proof} By \cite[Proposition A7]{IK} it is enough to prove that $f_1$ and $f_2$ define projectively equivalent 3-folds in $\PP^4$. 
Write $q_0(u,v)=\lambda _0p_0(u,v)+\lambda _1p_1(u,v)+\lambda _2p_2(u,v)$, $q_1=\mu _0p_0(u,v)+\mu _1p_1(u,v)+\mu _2p_2(u,v)$, $q_3=\rho _0p_0(u,v)+\rho _1p_1(u,v)+\rho _2p_2(u,v)$. We have 
$$\begin{array}{rcl} f_2 & = & q_0(u,v)x_0+q_1(u,v)x_1+q_2(u,v)x_2 \\
& = & (\lambda _0p_0(u,v)+\lambda _1p_1(u,v)+\lambda _2p_2(u,v))x_0+(\mu _0p_0(u,v)+\mu _1p_1(u,v)+\mu _2p_2(u,v))x_1+ \\
& & +(\rho _0p_0(u,v)+\rho _1p_1(u,v)+\rho _2p_2(u,v))x_2 \\ & = & (\lambda _0 x_0+\mu _0x_1+ \rho _0x_2)p_0(u,v)+ (\lambda _1 x_0+\mu _1x_1+ \rho _1x_2 ) p_1(u,v)+\\
& & +  (\lambda _2 x_0+\mu _2x_1+ \rho _2 x_2 ) p_2(u,v).
\end{array}
$$
Therefore, $f_1$ and $f_2$ define projectively equivalent hypersurfaces in $\PP^4$.
\end{proof}

\begin{remark} We fix  integers $d\ge 5$ and $2\le k \le \lfloor \frac{d}{2} \rfloor$.
We keep the notations introduced in Section 3. If $\rank M_k=3$, then $\rank \mathcal{A}_k \le 3$, $\rank \mathcal{B}_k\le 3$ and $\rank \mathcal{C}_k \le 3$.
\end{remark}

We will now explain  the geometrical meaning of the rank of the matrices $\mathcal{A}_k, \mathcal{B}_k, \mathcal{C}_k$ introduced in Section 3. To this end, we recall some basic facts about symmetric tensors in two variables. For more details see \cite{IK}, \cite{GO} and \cite{BGI}.

Let us fix an integer $t\geq 3$ and consider the vector space $K[u,v]_t$ of forms of degree $t$. Its elements can also be interpreted as symmetric tensors in two variables; by definition the Waring rank, or symmetric rank, of $p\in K[u,v]_t$ is the minimum integer $r$ such that there exist linear forms $l_1,\ldots, l_r\in K[u,v]_1$ such that $p=l_1^t+\cdots +l_r^t$. In particular, a symmetric tensor $p$ has Waring rank $1$ if $p=l^t$ for a suitable linear form $l$, i.e. $p$ is a pure power of degree $t$. 

In the projective space $\PP^t$, naturally identified with $\PP(K[u,v]_t)$,  the set of (equivalence classes of) forms of Waring rank $1$  is the image of the $t$-tuple Veronese embedding of $\PP^1$ in $\PP^t$, that is the rational normal curve  $C_t$ of degree $t$. We recall that, for any $r\geq 1$, the $r$-secant variety of $C_t$ is 
\[
\sigma_r(C_t)=\overline{\cup_{p_1,\ldots,p_r\in C_t}\langle p_1,\ldots,p_r\rangle}.
\]
Clearly $C_t=\sigma_1(C_t)\subset \sigma_2(C_t)\subset \cdots$, and a general element of $\sigma_r(C_t)\setminus \sigma_{r-1}(C_t)$ is a symmetric tensor of Waring rank $r$, but if $r>1$ $\sigma_r(C_t)$ contains also tensors of Waring rank $>r$.  The dimension of $\sigma_r(C_t)$ is $\min\{2r-1, t\}.$ Moreover, for any $r<\frac{t+1}{2}$,  $\sigma_{r-1}(C_t)$ is the singular locus of $\sigma_r(C_t)$ (see \cite[Proposition 1.2.2 and Corollary 1.2.3]{R}).

We recall also that the tangential surface of $C_t$, $TC_t$, is the closure of the union of the embedded tangent lines to $C_t$. The tangent line at the point $l_1^t\in C_t$ is the set of tensors that can be written in the form $l_1^{t-1}l_2$, with $l_2$ a linear form. Similarly the osculating $3$-fold of $C_t$, $T^2C_t$, is the closure of the union of the embedded osculating planes to $C_t$, and the osculating plane at $l_1^t$ is the set of tensors that can be written in the form $l_1^{t-2}m$, with $m$ a form of degree $2$. We are now ready to give the desired interpretation of the rank of the matrices introduced in Section 3. We state and prove  Proposition \ref{rank} for the form $p_0$ and the matrices $\mathcal{A}_k$; the analogous results hold true also for $p_1, p_2$, and their respectively catalecticant  matrices $\mathcal{B}_k$, $\mathcal{C}_k$.

\begin{proposition}\label{rank} We fix an integer $d\ge 5$ and we keep the notations introduced in Section 3. It holds: 
\begin{itemize}
    \item[(1)] If $\rank \mathcal{A}_k=1$ for some $2\le k \le \lfloor \frac{d+1}{2} \rfloor$ (and, hence, for all $k$), then $p_0=\ell ^{d-1}$ for some $\ell \in K[u,v]_1$.
    \item[(2)] If $\rank \mathcal{A}_k=2$ for some $3\le k \le \lfloor \frac{d+1}{2} \rfloor$ (and, hence, for all $k$), then either $p_0=\ell _1 ^{d-1}+\ell _2 ^{d-1}$ or $p_0=\ell _1^{d-2}\ell _2$ for some $\ell_1, \ell _2 \in K[u,v]_1$.
    \item[(3)] If $\rank \mathcal{A}_k=3$ for some $4\le k \le \lfloor \frac{d+1}{2} \rfloor$ (and, hence, for all $k$), then either $p_0=\ell _1 ^{d-1}+\ell _2 ^{d-1}+(\lambda \ell _1+\mu \ell _2)^{d-1}$ or $p_0=\ell _1^{d-1}+\ell _2^{d-2}(\lambda \ell _1+\mu \ell_2)$ for some $\ell_1, \ell _2 \in K[u,v]_1$ and $\lambda, \mu \in K^*$.
\end{itemize}
\end{proposition}
\begin{proof}
Let $r$ be any integer such that $r+1\leq k$. From \cite[Theorem 1.3]{GO},  it follows that all the minors of order $r+1$ of $\mathcal{A}_k$ vanish if and only if $[p_0]\in\sigma_r(C_{d-1})$. For $r=1$, this gives (1). For $r=2$, we get that if $\mathcal{A}_k$ has rank $2$, then  $p_0\in\sigma_2(C_{d-1})$. From \cite[Corollary 26]{BGI},  it follows that either $p_0$ has Waring rank $2$ or $p_0\in TC_{d-1}$; this proves (2). Similarly, for $r=3$, $\rank \mathcal{A}_k=3$ implies that $p_0\in \sigma_3(C_{d-1})$. So, either the Waring rank of $p_0$ is $3$, or $p_0$ belongs to the join of $C_{d-1}$ and its tangential surface (\cite[Corollary 26]{BGI}). This proves (3).
\end{proof}

\begin{theorem} \label{degree_d} The Artinian Gorenstein algebra $S/ \Ann _S(f)$ associated to a Perazzo 3-fold of degree $d\ge 5$ has $h$-vector: $(1, \ 5, \ 6, \ 6, \ \cdots \ 6, \ 6, \ 5, \ 1)$ if and only if, after a possible change of coordinates, one of the following cases holds:
\begin{itemize}
    \item[(i)]   $f(x_0,x_1,x_2,u,v)=u^{d-1}x_0+u^{d-2}vx_1+u^{d-3}v^2x_2+ au^d+bu^{d-1}v+cu^{d-2}v^2$ with $a,b,c\in K$, or
    \item[(ii)]  $f(x_0,x_1,x_2,u,v)=u^{d-1}x_0+u^{d-2}vx_1+v^{d-1}x_2+au^d+bu^{d-1}v+cv^{d}$ with $a,b,c\in K$, or
    \item[(iii)]  $f(x_0,x_1,x_2,u,v)=u^{d-1}x_0+(\lambda u+\mu v)^{d-1}x_1+v^{d-1}x_2+au^d+b(\lambda u+ \mu v)^d+cv^d$ with $\lambda, \mu\in K^*$ and $a,b,c\in K$.
\end{itemize}
\end{theorem}
\begin{proof}
As observed in Remark \ref{rank equal 3}, the  $h$-vector is minimal if and only if $\rank M_k=\rank N'_k=3$ for any $k$. A straigthforward computation shows that for any $f$ as in (i), (ii) or (iii) one has $\rank M_k=\rank N'_k=3$ for any $k$ and, therefore, $S/\Ann_S(f)$ has $h$-vector $(1, \ 5, \ 6, \ \cdots \ 6, \ 5, \ 1)$. To prove the converse, we first observe that if $\rank M_k=\rank N'_k=3$ for any $k$, then the ranks of $\mathcal{A}_k, \mathcal{B}_k, \mathcal{C}_k, \mathcal{G}_k$ are all bounded above by $3$. We analyze first the various possibilities for $p_0, p_1, p_2$.

(I) $d\geq 7$ and $\rank \mathcal{A}_k=\rank \mathcal{B}_k=\rank \mathcal{C}_k=3$ for $4\leq k\leq \lfloor\frac{d+1}{2}\rfloor $ and $p_0, p_1, p_2$ all have Waring rank $3$. We use  \cite[Corollary 1.2]{GO}:  the spaces of the columns of $\mathcal{A}_k, \mathcal{B}_k, \mathcal{C}_k$ coincide, so there exist linear forms $l_1, l_2, l_3$ and suitable constants such that 
$$p_0=\lambda_0l_1^{d-1}+\mu_0l_2^{d-1}+\nu_0l_3^{d-1}$$
$$p_1=\lambda_1l_1^{d-1}+\mu_1l_2^{d-1}+\nu_1l_3^{d-1}$$
$$p_2=\lambda_2l_1^{d-1}+\mu_2l_2^{d-1}+\nu_2l_3^{d-1}.$$
Since $p_0, p_1, p_2$ are linearly independent, the matrix
$\begin{pmatrix}
    \lambda_0 & \mu_0 & \nu_0\\
    \lambda_1 & \mu_1 & \nu_1\\
    \lambda_2 & \mu_2 & \nu 2
\end{pmatrix}$
is invertible, so $\langle p_0, p_1, p_2\rangle = \langle l_0^{d-1}, l_1^{d-1}, l_2^{d-1}\rangle$. In view of Lemma \ref{generators} in $f$ we can replace $p_0, p_1, p_2$ with $l_0^{d-1}, l_1^{d-1}, l_2^{d-1}$.

(II) $d\geq 7$ and $\rank \mathcal{A}_k=3$ for $4\leq k\leq \lfloor\frac{d+1}{2}\rfloor $, but $p_0$ has Waring rank strictly $>3$. So from Proposition \ref{rank} (3), $p_0$ is of the form $\ell _1^{d-1}+\ell _2^{d-2}(\alpha \ell _1+\beta \ell_2)$ for some $\ell_1, \ell _2 \in K[u,v]_1$ and $\alpha, \beta \in K^*$.   So up to the change of variables that sends $l_1$ into $u$, and $l_2$ into $v$, $p_0=u^{d-1}+\alpha uv^{d-2}+\beta v^{d-1}$. Then
\[ M_3=
\begin{pmatrix}
 1 & 0 & 0 & b_0 & b_1 & b_2 & c_0 & c_1 & c_2\\
 0 & 0 & 0 &     &&&&&\\
 \vdots & &&&&&&&\vdots\\
 0 & 0 & \alpha & b_{d-3} & b_{d-2} & b_{d-1} & c_{d-3} & c_{d-2} & c_{d-1} \\
 0 & \alpha & \beta & b_{d-2} & b_{d-1} & b_d & c_{d-2} & c_{d-1} & c_d
\end{pmatrix}.
\]
From $\rank M_3<4$ it follows $b_1=\cdots=b_{d-3}=c_1=\cdots=c_{d-3}=0$. Therefore
$$p_1=b_0u^{d-1}+b_{d-2}uv^{d-2}+b_{d-1}v^{d-1},\  p_2=c_0u^{d-1}+c_{d-2}uv^{d-2}+c_{d-1}v^{d-1},$$
and we can replace $p_0, p_1, p_2$ with $u^{d-1}, uv^{d-2}, v^{d-1}$.

(III) $\rank \mathcal{A}_k=2$ for $3\leq k\leq \lfloor\frac{d+1}{2}\rfloor $ and $p_0$ has Waring rank $2$, so it can be written $p_0=u^{d-1}+v^{d-1}$.
Then $M_3$ is as in case (II) with $\alpha=0$, $\beta=1$ 
and 
\[
M_2=
\begin{pmatrix}
 1 & 0 &b_0 & b_1 & c_0 & c_1\\
 0&0&b_1&b_2&c_1&c_2\\
 \vdots & \vdots & \vdots & \vdots & \vdots & \vdots\\
 0& 1 & b_{d-2} & b_{d-1}& c_{d-2} & c_{d-1}
\end{pmatrix}.
\]
From $\rank M_2<4$ we deduce that 
\[\rank \begin{pmatrix}
 b_1 & b_2\\
 \vdots & \vdots\\
 b_{d-3} & b_{d-2}
\end{pmatrix}<2, \ 
\rank \begin{pmatrix}
 c_1 & c_2\\
 \vdots & \vdots\\
 c_{d-3} & c_{d-2}
\end{pmatrix}<2, \ 
\rank \begin{pmatrix}
 b_1 & b_2 & \ldots & b_{d-2}\\
 c_1 & c_2 & \ldots & c_{d-2}
\end{pmatrix}<2.
\]
Therefore $$(b_1, \ldots, b_{d-2})=(\lambda^{d-3}, \lambda^{d-4}\mu, \ldots,\mu^{d-3}), \ (c_1, \ldots, c_{d-2})=(\sigma^{d-3}, \sigma^{d-4}\rho, \ldots,\rho^{d-3}),$$ for suitable $\lambda, \mu, \sigma, \rho\in K$. We get:
$$p_1=b_0u^{d-1}+uv((d-1)\lambda^{d-3}u^{d-3}+{\binom{d-1}{2}}\lambda^{d-4}\mu u^{d-4}v+\cdots)+b_{d-1}v^{d-1},$$
$$p_2=c_0u^{d-1}+uv((d-1)\sigma^{d-3}u^{d-3}+{\binom{d-1}{2}}\sigma^{d-4}\rho u^{d-4}v+\cdots)+c_{d-1}v^{d-1}.$$
We can also write
$$p_1=b_0u^{d-1}+uv\phi_{d-3}+b_{d-1}v^{d-1}, \ 
p_2=c_0u^{d-1}+kuv\phi_{d-3}+c_{d-1}v^{d-1}$$ where $\phi_{d-3}$ is a form of degree $d-3$ and $k\in K$, because $(b_1, \ldots,b_{d-2})$ and $(c_1, \ldots, c_{d-2})$ are proportional. We  can assume $b_0=c_0=0$ and we get $v^{d-1}\in \langle p_0, p_1, p_2\rangle$, hence $u^{d-1}, uv\phi_{d-3}\in \langle p_0, p_1, p_2\rangle$. Finally, adding to $uv\phi_{d-3}$ suitable multiples of $u^{d-1}, v^{d-1}$, we get $(\lambda u+\mu v)^{d-1}\in \langle p_0, p_1, p_2\rangle$.

(IV) $\rank \mathcal{A}_k=2$ for $3\leq k\leq \lfloor\frac{d+1}{2}\rfloor $ but $p_0$ has Waring rank $>2$, so up to a change of variables $p_0=u^{d-2}v$.
\[ M_2=
\begin{pmatrix}
0&1&b_0&b_1&c_0&c_1\\
1&0&b_1&b_2&c_1&c_2\\
\vdots & \vdots &&&&\\
0 & 0 & b_{d-2}& b_{d-1} & c_{d-2} & c_{d-1}
\end{pmatrix}
\]
has rank $3$, therefore
\[
\rank \begin{pmatrix}
b_2 & b_3 & c_2 & c_3\\
\vdots & \vdots & \vdots & \vdots\\
b_{d-2} & b_{d-1} & c_{d-2} & c_{d-1}
\end{pmatrix}<2,
\]
and arguing in a similar way to (III), we conclude that $\langle p_0, p_1, p_2\rangle$ is either of the form $\langle u^{d-1}, u^{d-2}v, (\lambda u+\mu v)^{d-1}\rangle $, or $\langle u^{d-1}, u^{d-2} v, u^{d-3}v^2\rangle.$

(V) $\rank \mathcal{A}_k=\rank \mathcal{B}_k=\rank \mathcal{C}_k=1$, then $p_0, p_1, p_2$ are all pure powers of degree $d-1$.

(VI)  Let $\pi$ be the $2$-plane generated by the polynomials $p_0, p_1, p_2$. If $d=5$, $\pi\subset \PP^4=\PP(K[u,v]_4)$. The tangential variety $TC_4$ has codimension $2$, so the intersection $\pi\cap TC_4\neq \emptyset$. If $\pi$ intersects $TC_4$ outside its singular locus $C_4$, up to a change of variables $u^3v\in\pi$ and we conclude as in (IV); otherwise, we are in the situation of (V). If $d=6$, $\pi\subset \PP^5=\PP(K[u,v]_5)$. Now $\sigma_2(C_5)$ has codimension $2$ and therefore $\pi\cap\sigma_2(C_5)\neq\emptyset$. Therefore we are either  in the situation of (III) or of (IV).

We have proved that for any $d\geq 5$, if $f$ defines a Perazzo $3$-fold and $S/\Ann_S(f)$ has minimal $h$-vector, then the polynomials $p_0, p_1, p_2$ are as in (i), or (ii), or (iii).

It remains to find out how we can choose the polynomial $g$ in each of the cases. From Proposition \ref{hilbert function} we deduce that the only condition that $g$ has to satisfy is $\Ann_S(p_0x_0+p_1x_1+p_2x_2)_3=\Ann_S(f-g)_3=\Ann_S(f)_3$. In other words, we impose that $g$ is annihilated by a system of generators of $\Ann_S(f)_3$.

(i) If $f(x_0,x_1,x_2,u,v)=u^{d-1}x_0+u^{d-2}vx_1+u^{d-3}v^2x_2+g$, we have that $\Ann_S(f)_3=\langle V^3\rangle$ and so $g=g_0u^d+g_1u^{d-1}v+g_2u^{d-2}v^2$. 

(ii) If $f(x_0,x_1,x_2,u,v)=u^{d-1}x_0+u^{d-2}vx_1+v^{d-1}x_2+g$, we have that $\Ann_S(f)_3=\langle UV^2\rangle$. This gives that $\sum_{i=2}^{d-1}g_i\binom{d}{i}(k-i)i(i-1)u^{k-i-1}v^{i-2}=0$, so $g_2=\dots=g_{d-1}=0$. Thus we get $g=g_0u^d+g_1u^{d-1}v+g_dv^d$.

(iii) If $f(x_0,x_1,x_2,u,v)=u^{d-1}x_0+(\lambda u+\mu v)^{d-1}x_1+v^{d-1}x_2+g$, we have that $\Ann_S(f)_3=\langle \mu U^2V-\lambda UV^2\rangle$. Then we have the condition $$\sum_{i=2}^{d-1}\binom{k-3}{i-1}(\mu g_i-\lambda g_{i+1})u^{d-k-2}v^{i-1}=0 \iff \mu g_i-\lambda g_{i+1}=0,\quad i=1,\dots,k-2.$$
So we can collect $g_1$ and complete the $d$-th power to obtain $g=au^d+b(\lambda u+\mu v)^d+cv^d$.
\end{proof}

\begin{remark}
As we noticed in Lemma \ref{generators}, the Hilbert function of the algebra $S/\Ann_S(f)$ depends only on the plane $\pi=\langle p_0, p_1, p_2\rangle\subset\PP^{d-1}$ and not on the choice of the three generators. In Theorem \ref{degree_d} we have proved that the $h$-vector is minimal if and only if the plane $\pi$ is in one of the following positions: it is an osculating plane to the rational normal curve $C_{d-1}$ (case (i)), or it contains the tangent line to $C_{d-1}$ at a point and meets $C_{d-1}$ also at a second point (case (ii)), or it intersects $C_{d-1}$ at three distinct points (case (iii)).
\end{remark}

\begin{remark}
    In Theorem \ref{degree_d}, we have obtained a complete characterization of the polynomials $f$ such that $A=S/\Ann_S(f)$ has minimum $h$-vector for any $d\geq 5$. This allows us to conclude with a direct verification  the proof of Theorem \ref{main2}, proving the WLP of these algebras in the cases $5\leq d\leq 7$.
\end{remark}

\begin{corollary}\label{quintic}
The Artinian Gorenstein algebra $S/ \Ann _S(f)$ associated to a Perazzo 3-fold of degree 5 has the WLP if and only if, after a possible change of coordinates, one of the following cases holds:
\begin{itemize}
    \item[(i)]   $f(x_0,x_1,x_2,u,v)=u^4x_0+u^3vx_1+u^2v^2x_2+au^5+bu^4v+cu^3v^2\in R_5$ with $a,b,c\in K$, or
    \item[(ii)]  $f(x_0,x_1,x_2,u,v)=u^4x_0+u^3vx_1+v^4x_2+au^5+bu^4v+cv^5\in R_5$ with $a,b,c\in K$, or
    \item[(iii)]  $f(x_0,x_1,x_2,u,v)= u^4x_0+(\lambda u+ \mu v)^4x_1+v^4x_2+au^5+b(\lambda u+ \mu v)^5+cv^5\in R_5$ with $\lambda, \mu\in K^*$ and $a,b,c\in K$.
\end{itemize}
\end{corollary}
\begin{proof}
It follows from Theorems \ref{degree_d}, \ref{main1}, and \ref{main2}.
\end{proof}

Note that, as consequence of the results of Gordan-Noether, Corollary \ref{quintic} gives also a complete classification of threefolds in $\PP^4$ with vanishing hessian.

\vskip 4mm
\section{Final comments}
In this last section, we give a short geometrical description of the hypersurfaces of   Theorem \ref{degree_d} when $a=b=c=0$.

Case (i) corresponds to the union of the classic cubic Perazzo 3-fold in $\PP^4$ of equation: $u^2x_0+uvx_1+v^2x_2=0$ with the non-reduced  hyperplane of equation: $u^{d-3}=0$. To describe the other two hypersurfaces, we first recall some known geometric properties of hypersurfaces with vanishing hessian. Let $X=V(f)\subset\PP^N$ be such a hypersurface with $\hess_f=0$.

We denote by 
$$\nabla _{f}:\PP^N \dashrightarrow (\PP^N)^* $$

\noindent its {\em polar map} defined  by
$$ \nabla _{f}(p)=\left(\frac{\partial f}{\partial x_0}(p),\frac{\partial f}{\partial x_1}(p), \ldots, \frac{\partial f}{\partial x_N}(p)\right),$$
and by 
$$\gamma :X \dashrightarrow (\PP^N)^*
$$
the restriction of $\nabla _{f}$ to $X$, i.e. the {\em  Gauss map} of $X$, associating to each smooth point of $X$ its embedded tangent space. 
The image of $\gamma$ is the dual variety $X^*$ of $X$. Let $Z=\overline{\nabla _f(\PP^N)}$ be the  closure of the   image of the polar map. Then $X^*\subsetneq Z \subsetneq (\PP^N)^*$ (\cite[Corollary 7.2.8]{R}). 
Moreover, if $N=4$, $Z$ is a cone with vertex a line over an irreducible plane curve, and its dual $Z^*$ is a rational plane curve in $\PP^4$, naturally identified with the bidual space $(\PP^4)^*{^*}$ (\cite[Lemma 7.4.13]{R}). 

Let  $X$ be a Perazzo hypersurface of degree $d$ in $\PP^4$ of equation (\ref{perazzo form}). $X$ contains the line $L:x_0=x_1=x_2=0$ and the plane $\Pi: u=v=0$. From \cite[Sections 7.3 and 7.4]{R}, it follows that $\Pi$ is the singular locus of $X$ with  multiplicity $d-1$; moreover, $X^*$ is a scroll surface of degree $d$, having the line $\Pi^*$ as directrix. In particular, $\Pi^*$ is also the vertex of $Z$, and the general plane ruling of the cone $Z$ meets $X^*$ along a line of the scroll. The curve $Z^*$ is contained in $\Pi$ and the hyperplanes containing $\Pi$ cut on $X$, outside $\Pi$, a $1$-dimensional family $\Sigma$ of planes: they are all tangent to  $Z^*$ and meet $L$. If $p$ is general in $X$, then the fibre of the Gauss map $\gamma^{-1}(\gamma(p))$ is the line $\langle p, p'\rangle$ where $p'$  is the tangency point to $Z^*$ of the plane of the  family $\Sigma$ passing through $p$.

We now see how this picture specializes if we consider the reduced, irreducible Perazzo 3-fold   $X_1\subset\PP^4$ of equation 
 $$ f_1(x_0,x_1,x_2,u,v)=u^{d-1}x_0+u^{d-2}vx_1+v^{d-1}x_2,$$ 
 case (ii) in Theorem \ref{degree_d}.
 We use coordinates $z_0,\ldots,z_4$ in $(\PP^4)^*$. The equation of $Z$, which expresses the algebraic dependence of $p_0, p_1,p_2$, is $z_1^{d-1}-z_0^{d-2}z_2=0$; the one of $Z^*$ is $(d-1)^{d-1}x_0^{d-2}x_2+(d-2)^{d-2}x_1^{d-1}=0$. They both represent rational curves of degree $d-1$ with a singular point of multiplicity $d-2$ with only one tangent line.
 
 In case (iii) we have
 $$ f_2(x_0,x_1,x_2,u,v)=u^{d-1}x_0+(\lambda u+\mu v)^{d-1}x_1+v^{d-1}x_2 \text{ with } \lambda, \mu \in K^*.$$
 
 For low values of $d$ we have checked with the help of Macaulay2 (\cite{M2}) that $Z$ is a cone over a rational  curve of degree $d-1$ with $\frac{(d-2)(d-3)}{2}$ distinct nodes. Its dual $Z^*$ results to be  a rational curve of degree $2d-4$. If $d=5$, then $Z^*$ has degree $6$ and it has $3$ cuspidal points of multiplicity $3$ at the fundamental points $[1,0,0], [0,1,0], [0,0,1]$ and one node; if $d=6$, then $Z^*$ has cuspidal points of multiplicity $4$ at the fundamental points and $3$ nodes; if $d=7$, then $Z^*$ has cuspidal points of multiplicity $5$ at the fundamental points and $6$ nodes.

\end{document}